
\frenchspacing

\parindent15pt

\abovedisplayskip4pt plus2pt
\belowdisplayskip4pt plus2pt
\abovedisplayshortskip2pt plus2pt
\belowdisplayshortskip2pt plus2pt

\font\twbf=cmbx10 at12pt
 at12pt
 at12pt

\font\sc=cmcsc10

\font\ninerm=cmr9
\font\nineit=cmti9
\font\ninesy=cmsy9
\font\ninei=cmmi9
\font\ninebf=cmbx9

\font\sevenrm=cmr7

\font\seveni=cmmi7
\font\sevensy=cmsy7

\font\fivenrm=cmr5
\font\fiveni=cmmi5
\font\fivensy=cmsy5

\def\nine{%
\textfont0=\ninerm \scriptfont0=\sevenrm \scriptscriptfont0=\fivenrm
\textfont1=\ninei \scriptfont1=\seveni \scriptscriptfont1=\fiveni
\textfont2=\ninesy \scriptfont2=\sevensy \scriptscriptfont2=\fivensy
\textfont3=\tenex \scriptfont3=\tenex \scriptscriptfont3=\tenex
\def\rm{\fam0\ninerm}%
\textfont\itfam=\nineit
\def\it{\fam\itfam\nineit}%
\textfont\bffam=\ninebf
\def\bf{\fam\bffam\ninebf}%
\normalbaselineskip=11pt
\setbox\strutbox=\hbox{\vrule height8pt depth3pt width0pt}%
\normalbaselines\rm}

\hsize30cc
\vsize44cc
\nopagenumbers

\def\luz#1{\luzno#1?}
\def\luzno#1{\ifx#1?\let\next=\relax\yyy
\else \let\next=\luzno#1\xxx\fi\next}
\def\sp#1{\def\xxx{\kern1.7pt}\def\yyy{\kern-1.7pt}\luz{#1}}
\def\spa#1{\def\xxx{\kern1pt}\def\yyy{\kern-1pt}\luz{#1}}

\newcount\beg
\newbox\aabox
\newbox\atbox
\newbox\fpbox
\def\abbrevauthors#1{\setbox\aabox=\hbox{\sevenrm\uppercase{#1}}}
\def\abbrevtitle#1{\setbox\atbox=\hbox{\sevenrm\uppercase{#1}}}
\long\def\pag{\beg=\pageno
\def\leftheadline{\noindent\rlap{\nine\folio}\hfil\copy\aabox\hfil}
\def\rightheadline{\noindent\hfill\copy\atbox\hfill\llap{\nine\folio}}
\def\phead{\setbox\fpbox=\hbox{\sevenrm
************************************************}%
\noindent\vbox{\sevenrm\baselineskip9pt\hsize\wd\fpbox%
\centerline{                                             }

\centerline{                                        }

\centerline{                                 }

\centerline{                            }

\centerline{                  }}\hfill}
\footline{\ifnum\beg=\pageno \hfill\nine[\folio]\hfill\fi}
\headline{\ifnum\beg=\pageno\phead
\else
\ifodd\pageno\rightheadline \else \leftheadline \fi
\fi}}

\newbox\tbox
\newbox\aubox
\newbox\adbox
\newbox\mathbox

\def\title#1{\setbox\tbox=\hbox{\let\\=\cr
\baselineskip14pt\vbox{\twbf\tabskip 0pt plus15cc
\halign to\hsize{\hfil\ignorespaces \uppercase{##}\hfil\cr#1\cr}}}}

\newbox\abbox
\setbox\abbox=\vbox{\vglue18pt}

\def\author#1{\setbox\aubox=\hbox{\let\\=\cr
\nine\baselineskip12pt\vbox{\tabskip 0pt plus15cc
\halign to\hsize{\hfil\ignorespaces \uppercase{\spa{##}}\hfil\cr#1\cr}}}%
\global\setbox\abbox=\vbox{\unvbox\abbox\box\aubox\vskip8pt}}

\def\address#1{\setbox\adbox=\hbox{\let\\=\cr
\nine\baselineskip12pt\vbox{\it\tabskip 0pt plus15cc
\halign to\hsize{\hfil\ignorespaces {##}\hfil\cr#1\cr}}}%
\global\setbox\abbox=\vbox{\unvbox\abbox\box\adbox\vskip16pt}}

\def\mathclass#1{\setbox\mathbox=\hbox{\footnote{}{1991 {\it Mathematics Subject
Classification}\/: #1}}}

\long\def\maketitlebcp{\pag\unhbox\mathbox
\vglue7cc
\box\tbox
\box\abbox
\vskip8pt}

\long\def\abstract#1{{\nine{\bf Abstract.} 
#1

}}

\def\section#1{\vskip-\lastskip\vskip12pt plus2pt minus2pt
{\bf #1}}

\long\def\th#1#2#3{\vskip-\lastskip\vskip4pt plus2pt
{\sc #1} #2\hskip-\lastskip\ {\it #3}\vskip-\lastskip\vskip4pt plus2pt}

\long\def\remar#1#2{\vskip-\lastskip\vskip4pt plus2pt
\sp{#1} #2\vskip-\lastskip\vskip4pt plus2pt}

\def\Proof{\vskip-\lastskip\vskip4pt plus2pt
\sp{Proo{f.}\ }\ignorespaces}

\def\endproof{\nobreak\kern5pt\nobreak\vrule height4pt width4pt depth0pt
\vskip4pt plus2pt}

\newbox\refbox
\newdimen\refwidth
\long\def\references#1#2{{\nine
\setbox\refbox=\hbox{\nine[#1]}\refwidth\wd\refbox\advance\refwidth by 12pt%
\def\textindent##1{\indent\llap{##1\hskip12pt}\ignorespaces}
\vskip24pt plus4pt minus4pt
\centerline{\bf References}
\vskip12pt plus2pt minus2pt
\parindent=\refwidth
#2

}}

\def\footnoterule{\kern -3pt \hrule width 4cc \kern 2.6pt}

\catcode`@=11
\def\vfootnote#1%
{\insert\footins\bgroup\nine\interlinepenalty\interfootnotelinepenalty%
\splittopskip\ht\strutbox\splitmaxdepth\dp\strutbox\floatingpenalty\@MM%
\leftskip\z@skip\rightskip\z@skip\spaceskip\z@skip\xspaceskip\z@skip%
\textindent{#1}\footstrut\futurelet\next\fo@t}
\catcode`@=12


\mathclass{Primary 57N10; Secondary 57M20.}

\abbrevauthors{M. Sokolov}
\abbrevtitle{Summands of the Turaev--Viro invariant}

\title{On the absolute value of the SO(3)--invariant\\
and other summands\\
of the Turaev--Viro invariant}

\author{Maxim\ Sokolov}
\address{Chelyabinsk State University\\
ul. Br. Kashirinykh, 129\\
454136 Chelyabinsk, Russia\\
E-mail: sok@cgu.chel.su}

\maketitlebcp

\footnote{}{The author is partially  supported by INTAS, grant 94201.}

\footnote{}{In this version we added a definition
of the invariant $TV_1$ on a triangulation,  Lemma 3 and
Remark 5 in section 6, some remarks about the tables.
We also corrected the definition
of an odd edge in section 6.}


\section{1. Introduction.} It was proved in [S1] and [S2] that each
Turaev--Viro invariant $TV(M)_q$ for a 3-manifold $M\/$ is a sum of three
invariants $TV_0(M)_q, TV_1(M)_q$, and $TV_2(M)_q$ (for definition of the
Turaev--Viro invariants, see [TV]). It  follows from the Turaev--Walker
theorem (see [T1], [W]) that if $q^2$ is a primitive root of unity of an odd
degree then, up to normalization,
$TV_0(M)_q$ coincides with  the square of the modulus of the  so--called
$SO(3)$-invariant $\tau_e(M)$ defined in [T2].
For  a connection between $SO(3)$-invariants and
the Reshetikhin--Turaev invariants, see [KM] and [BHMV].

With a help of suitable normalizations we make the numbers
$(TV_0(M)_q+TV_2(M)_q)$, $TV_0(M)_q$, and $TV_1(M)_q$
to be invariant under removing of 3-balls.
That allows us to
define these three invariants on a triangulation of a closed 3-manifold $M\/$.

It is natural  to relate the invariants $TV_N(M),\quad N=0,1,2,$
 to the Turaev--Viro
invariants. Here we show that for every 3-manifold $M\/$ the following
 holds:
$$
TV_0(M)_q+TV_2(M)_q={1\over 2} (TV(M)_q+ TV(M)_{-q}),
$$
$$
TV_1(M)_q={1\over 2}(TV(M)_q-TV(M)_{-q}).
$$

At the end of the paper we present a few tables. There are a lot of
numerical
tables of the Turaev--Viro and Reshetikhin--Turaev invariants (see,
for instance, [KL1],
[KL2], [N], [S2]). An advantage of our tables is that the values
of the invariants are presented
as polynomials on $q\/$ with integer coefficients.

I would like to express my thanks to my supervisor Prof. S.~V.~Matveev,
to Prof. R.~Zh.~Aleev and
to Prof. C.~Blanchet for many useful conversations. An interesting question
of Prof. D.~Yetter was a reason for appearing of  section 5.
A. Macovetsky helped me to prove lemma 3 in section 6. I am greatly indebted
to Prof. G. Masbaum who helped to find some deficiencies in the first version
of the paper.
My special
thanks to the Organizing committee of the mini-semester on Knot theory
in Warsaw, July--August 1995.

\section{2. Simple polyhedra and their local moves.}
A 2-dimensional polyhedron $X$ is called {\it simple} if the link of any
point of $X$ is homeomorphic to one of the following polyhedra:
(1) a circle, (2) a circle with two radii,
(3) a circle with three radii, (4) the segment $[0,1]$, (5) a wedge of
three segments with a common endpoint.

The set of points of a simple polyhedron $X$  having  links of types
(4) or (5) is called a {\it boundary} of $X$ and denoted by $\partial X$.
The points with links of   type (3) are called  {\it vertices} of $X$.
By an {\it edge} of $X$ we mean a connected component of the set of points
having the links of  type (2).

Simple polyhedra are also called fake surfaces. This class of polyhedra
generates the class of  special
polyhedra. Recall that a simple polyhedron $X$ is called  {\it special}
 if $\partial X=\emptyset$ and each 2-component of $X$ is a 2-cell.

A simple polyhedron $X$ with $\partial X = \emptyset$
is called a {\it simple spine} of a compact 3-manifold
$M$ with $\partial M \not =\emptyset$ if there exists an embedding $i\colon
X\rightarrow M$ such that $M \searrow i(X)$, i.e. $M$ collapses onto $i(X)$.
In the case of a closed $M$, a polyhedron $X$ is called a {\it simple spine} of
$M$ if it is a simple spine of $M$ with an open 3-ball removed. A simple
spine is called  {\it special} if it is a special polyhedron. It is
known that every compact connected 3-manifold has a special spine (see [Ca],
[M]).

Let as describe now  special polyhedra-with-boundary $P_1,\ldots ,P_4$.
Let $P_1$ be the polyhedron obtained from a disk $D^2\/$ by attaching two
semidisks along two parallel chords, $h_1$ and $h_2$ of $D^2$. The polyhedron
$P_2$ is obtained from $D^2\/$ by attaching a semidisk along $h_2$ and
the second one along a simple curve $l\/$ in $D^2\/$ that has the same endpoints
as $h_1$ and intersects $h_2$ transversally in exactly two points.
Let $R=R_1\cup R_2\cup R_3$ be a triod consisting of three radii of the
disk $D^2$. The polyhedron $P_3$ is obtained from the polyhedron
$(D^2\times \{ 0\})\cup (R\times I)$ by attaching a semidisk along a chord
$h_1\subset D^2\/$ that intersects the radius $R_1$ in just one interior
point. The polyhedron $P_4$ is obtained from
$(D^2\times \{ 0\})\cup (R\times I)$ by attaching a semidisk along
a simple curve that has the same endpoints as $h_1$ and intersects the triod
$R\/$ in exactly two points, on $R_2\/$ and $R_3$.

By {\it ${\cal L}$-move\/} on  simple polyhedra we
mean a replacement of a
fragment  homeomorphic to $P_1$ by  $P_2$.
By {\it ${\cal M}$-move\/} on simple polyhedra we mean a replacement
of a
fragment homeomorphic to $P_3$ by  $P_4$
(for details, see  [M],[P]).

Let a circle $c\/$ bounds a 2-disk in a 2-component of a special polyhedron
$X$. By {\it ${\cal B}$--move\/} we mean an attaching of additional 2-disk
to $X\/$ along $c$ (for details, see  [TV]).

It is proved in [M] that any two special spines of a 3-manifold can be
transformed one to another by a sequence of the moves ${\cal M}^{\pm 1}$ and
${\cal L}^{\pm 1}$. Note, that applying $\cal L$ several times, one can
transform any simple spine into a special one. So the theorem of
S.~V.~Matveev is true for simple spines too (see [P] also).

The $\cal B$--move on a simple spine of a 3-manifold $M\/$ corresponds to
removing of one 3-ball from $M\/$.

\section{3. The Turaev--Viro invariants.}
Throughout the paper, let us fix $r\ge 3$ and a root of unity $q$ of degree
$2r$ such that $q^2$ is a primitive root  of degree $r$.

In this section we recall how V.~G.~Turaev and O.~Y.~Viro define their
invariants on a simple polyhedron $X$ (cf. [TV]). Let $v_1,\ldots , v_d$
be the vertices of $X$, let $e_1,\ldots , e_f$ be the edges of $\partial X$
and let $\Gamma_1,\ldots , \Gamma_b$ be the 2-components of $X$.

By a {\it coloring\/}  of  $X\/$  we
mean            an             arbitrary             mapping
$$
\varphi\colon\{\Gamma_1,\ldots,\Gamma_b\}\to
{\bf Z}_{r-1}=\{ 0,1,\ldots ,r-2\}.
$$
A triple $(i,j,k)\in {\bf Z}^3_{r-1}$ will be called {\it admissible\/} if
$$
2r-4 \geq i+j+k
\equiv 0\pmod2,
$$
$$
\vert i-j \vert\leq
k \leq i+j.
$$
A coloring $\varphi\/$ is called {\it admissible\/} if for any edge $E\/$
of $X-\partial X\/$ the colors of the 2-components incident to $E\/$ form
an admissible triple.
Let us denote the set of admissible triples by $adm$ and
the set of admissible colorings of $X\/$ by  $Adm(X).$

By a {\it coloring\/} of a regular graph $G\/$ we shall mean any
mapping of the set of its  edges to ${\bf Z}_{r-1}$.
Let us denote the set of colorings of $X\/$ by $Col(X)$. Any coloring
$\varphi\/$ of a simple polyhedron $X\/$ induces in a natural way a
coloring $\partial \varphi\/$ of its boundary $\partial X\/$: an edge of
$\partial X\/$ takes the color of the 2-component of $X\/$ in whose boundary
this edge is contained.

Let $\Gamma_i,\Gamma_j,\Gamma_k$ be 2-components incident
to an edge  $E\/$  of  $X\/$
and   let   $\varphi\in Adm(X)$.   We   shall   say   that an unordered triple
$\{\varphi(\Gamma_i),\varphi(\Gamma_j),\varphi(\Gamma_k)\}$ is
a {\it color of the edge}  $E$.  There  are  six {\it wings\/}
incident to any vertex $v$ of  a  simple polyhedron.
Suppose   they  receive  under  $\varphi\/$  the
values $i,j,k,l,m,n\in{\bf Z}_{r-1}$. A 6-tuple
$\pmatrix{i&j&k\cr
          l&m&n\cr}$  is  called  a
{\it color of the vertex} $v\/$ if  $\{ i,j,k\}$  is a color  of  some  edge
incident to $v\/$ and $(i,l)$, $(j,m)$, $(k,n)\/$ are the   colors
of opposite 2-components incident to $v$.

     For an integer $n>0$ set
$$[n]_q    ={q^n-q^{-n} \over q-q^{-1}},$$
$$[n]_q   ! = [n]_q[n-1]_q\ldots [2]_q[1]_q.$$
Set also $[0]_q   =[0]_q   ! =1$. For a color $\{ i,j,k\}$ of  an  edge
set
$$
\Delta_q (i,j,k)=\left([\underline i +\underline j -\underline k]_q!
[\underline i+\underline k-\underline j]_q!
[\underline j+\underline k-\underline i]_q!\over
[\underline i+\underline j+\underline k+1]_q!\right)^{1/2}$$
where $\underline i=i/2$. Note that the expression in the round  brackets
presents a real number. By the square root $x^{1/2}$     of  a  real
number $x\/$ we mean the positive root of $|x|$ multiplied by $\sqrt{-1}$
if $x<0$.

     Let $\pmatrix{i&j&k\cr
          l&m&n\cr}$ be a color of some vertex $v$. A {\it symbol\/}  of
$v\/$ is defined by the following formula
$$
\vert T^{\varphi}_v\vert_q =\left\vert\matrix{i&j&k\cr
                                              l&m&n\cr}\right\vert_q
                                          =(\sqrt{-1})^{-(i+j+k+l+m+n)}
\Delta_q (i,j,k)\Delta_q (i,m,n)\times
$$
$$
\hfill{}\times\Delta_q (j,l,n)\Delta_q(k,l,m)
\left [\matrix{i&j&k\cr
               l&m&n\cr}\right ]_q,
$$
where
$$
\left[\matrix{i&j&k\cr
              l&m&n\cr}\right]_q =
\sum_z(-1)^z [z+1]_q!
\{ [z-\underline i-\underline j-\underline k]_q!
[z-\underline i-\underline m-\underline n]_q!
[z-\underline j-\underline l-\underline n]_q!
[z-\underline k-\underline l-\underline m]_q!\times
$$
$$
\hfill{}\times[\underline i+\underline j+\underline l+\underline m-z]_q!
[\underline i+\underline k+\underline l+\underline n-z]_q!
[\underline j+\underline k+\underline m+\underline n-z]_q!\}^{-1}.\quad
$$
Here $z\/$ runs through the non-negative  integers  such
that  all   expressions   in   the   square   brackets   are
non-negative. For $i\in {\bf Z}_{r-1}$ put
$$
w_{i,q} = (\sqrt {-1})^i [i+1]_q^{1/2}.$$
For $\varphi\in Adm(X)\/$ put
$$
|X,\varphi|_q = \prod_{i=1}^b w^{2\chi(\Gamma_i)}_{\varphi (\Gamma_i),q}
\prod_{s=1}^f w^{\chi(e_s)}_{\partial\varphi (e_s),q}
\prod_{j=1}^d|T^{\varphi}_{v_j}|_q,
$$
where $\chi$ is the Euler characteristic (the 2-components of $X\/$ and
the edges of $\partial X\/$ are thought to be open,
so if $e_s$ is homeomorphic to
${\bf R}$ then $\chi(e_s)=-1$ and if $e_s$ is homeomorphic to $S^1\/$ then
$\chi(e_s)=0$).

The Turaev-Viro invariant for the simple polyhedron $X\/$ is given by
$$
TV(X)_q=\sum_{\varphi\in Adm(X)}|X,\varphi|_q.
$$
It is proved in [TV] that $TV(X)_q$ is invariant under moves
${\cal L}^{\pm 1}$ and ${\cal M}^{\pm 1}$. It follows from Matveev's
theorem that if $X$ is a simple spine of a 3-manifold $M$ then
$TV(M)_q=TV(X)_q$ is a topological invariant of $M$.

Note that in [TV] a different normalization is used. The original
Turaev-Viro invariant
is given by the formula
$$
TV^*(X)_q=\omega^{-2\chi(X)+\chi(\partial X)}TV(X)_q ,
$$
where $\omega=\sqrt{2r}/|q-q^{-1}|$. It is proved in [TV] that $TV^*(X)_q$
is  invariant under ${\cal B}^{\pm 1}$ also.

\remar{Remark\ {1.}\ }{It is easily seen that if $q$ is a primitive root of
unity of degree $2r$ and
$\partial X=\emptyset$ then the numbers $|X,\varphi|_q$, and therefore the
numbers $TV(X)_q$ and $TV^*(X)_q$, lie in ${\bf Q}(q)$.}

\section{4. The summand--invariants.}
The set of 2-components of $X\/$ that receive odd colors under a coloring
$\varphi\in Adm(X)$ forms a surface embedded in $X\/$. We  denote this
surface by $S(\varphi)$. Note that $\partial S(\varphi)\subseteq\partial X$.

Present the set $Adm(X)$ as a disjoint union of subsets $Adm_0(X), Adm_1(X)$
and $Adm_2(X)$, where
\item{0)} $\varphi\in Adm_0 (X)\Leftrightarrow
(\varphi\in Adm(X))$ \& $(S(\varphi)=
    \emptyset)$;
\item{1)} $\varphi\in Adm_1( X)\Leftrightarrow
(\varphi\in Adm(X))$ \& $(\chi (S(\varphi))
    \equiv 1\pmod 2)$;
\item{2)} $\varphi\in Adm_2( X)\Leftrightarrow (\varphi\in Adm(X))$ \&
$(S(\varphi)\not=\emptyset)$ \&
$(\chi (S(\varphi))\equiv 0\pmod 2)$.

For any coloring $\alpha$ of $\partial X\/$ and $N\in\{ 0,1,2\}$ put
$$
\Omega_N(X,\alpha)_q=\sum_{\scriptstyle \varphi\in Adm_N(X)\atop
\scriptstyle \partial\varphi=\alpha} |X,\varphi|_q .
$$
If $\{\varphi\in Adm_N(X)\colon\partial\varphi =\alpha\}=\emptyset$, then
$\Omega_N(X,\alpha)_q=0$. Put also
$$
TV_N(X)_q=\sum_{\alpha\in Col(\partial X)} \Omega_N(X,\alpha)_q,
$$
where sum is taken over all colorings of $\partial X\/$.

\remar{Remark\ {2.}\  }{$TV(X)_q=TV_0(X)_q+TV_1(X)_q+TV_2(X)_q.$}

\remar{Remark\ {3.}\  }{If $q$ is a primitive root of unity of degree $2r$
then for a simple polyhedron $X\/$ with $\partial X
=\emptyset$ we have $TV_N(X)_q\in {\bf Q}(q)$, for any $N\in \{ 0,1,2 \}$
(see remark 1).}

\th{Lemma}{1.}{Let a simple polyhedron $X\/$ be the union of simple polyhedra
$Y\/$ and $Z\/$ and let each connected component of $T=Y\cap Z\/$ be a
connected component of both $\partial Y\/$ and $\partial Z\/$. Then for
any coloring $\beta\/$ of $\partial X\/$ we have
$$
\Omega_0(X,\beta)_q=\>\sum_{\alpha\in Col(T)}\quad
\Omega_0(Y,\alpha\cup(\beta|_{Y\cap\partial X}))_q\,
\Omega_0(Z,\alpha\cup(\beta|_{Z\cap\partial X}))_q,
$$

$$
\Omega_1(X,\beta)_q=\>\sum_{\scriptstyle \alpha\in Col(T) \atop
\scriptstyle K+L\equiv 1 (2)}
\Omega_K(Y,\alpha\cup(\beta|_{Y\cap\partial X}))_q\,
\Omega_L(Z,\alpha\cup(\beta|_{Z\cap\partial X}))_q,
$$

$$
\Omega_2(X,\beta)_q=\sum_{\scriptstyle \alpha\in Col(T) \atop
\scriptstyle K+L=2 \, {or}\, 4}
\Omega_K(Y,\alpha\cup(\beta|_{Y\cap\partial X}))_q\,
\Omega_L(Z,\alpha\cup(\beta|_{Z\cap\partial X}))_q.
$$
}

\Proof This follows from the equalities
$$
\vert X,\varphi\vert_q =\vert Y,(\varphi|_Y)\vert_q \cdot
\vert Z,(\varphi|_Z)\vert_q,
$$
 where $\varphi\in Adm(X)$
(see Lemma 4.2.A in [TV]), and $\chi (X)=\chi (Y)+\chi (Z)$.
\vskip4pt plus2pt

\th{Theorem}{1.}{ Let $X$ be a simple 2-polyhedron and $\alpha$ be a coloring
of $\partial X$. Then for any $N\in \{ 0,1,2 \}$ the number
$\Omega_N(X,\alpha)_q$ is invariant under ${\cal L}^{\pm 1}$ and
${\cal M}^{\pm 1}$.}

\Proof Let us show that the number $\Omega_N(X,\alpha)_q$ is invariant under
$\cal L$. The case of $\cal M$-move is  similar. By lemma 1 it is
sufficient to prove that $\Omega_N(P_1,\gamma)_q=\Omega_N(P_2,\gamma)_q$
for any $N\in\{ 0,1,2\}$,
where $P_1$ and $P_2$ are the polyhedra from the definition of the
${\cal L}$-move and $\gamma\/$
is a coloring of the graph $\partial P_1=\partial P_2$. It is easy to check
that for any $\gamma\/$ there is a unique $K\in\{ 0,1,2\}$ such that
$\{ \varphi\in Adm(P_i)\colon\partial\varphi=\gamma\}\subset Adm_K(P_i)$,
for $i=1, 2$. Therefore
$$
\Omega_N(P_1,\gamma)_q=
\sum_{\scriptstyle \varphi\in Adm(P_1)\atop
\scriptstyle \partial\varphi=\gamma}|P_1,\varphi|_q\quad {\rm and}\quad
\Omega_N(P_2,\gamma)_q=
\sum_{\scriptstyle \psi\in Adm(P_2)\atop
\scriptstyle \partial\psi=\gamma}|P_2,\psi|_q
$$
if $N=K$, and $\Omega_N(P_1,\gamma)_q=\Omega_N(P_2,\gamma)_q= 0$ if $N\not= K$.
It is proved in Lemma 4.4.A of [TV] that the sums are equal.
\vskip4pt plus2pt

\th{Corollary}{1.}{Let $X\/$ be a simple spine of a 3-manifold $M$.
Then $TV_N(M)_q=TV_N(X)_q$ is an invariant of $M$ for any $N\in\{ 0,1,2\}$.}
\vskip4pt plus2pt

\section{5. The summand-invariants and a triangulation.}
The summand invariants are not invariants under $\cal B$--move. This
prevents us from defining these invariants on a triangulation of a 3-manifold.
 Here
we modify the invariants $TV_0$, $TV_1$, and $TV_0+TV_2$ to make
them invariant
under removing of 3-balls.

Put $\omega_0=\sqrt{r}/|q-q^{-1}|$ and
$\omega =\sqrt{2r}/|q-q^{-1}|$. Let $X\/$ be a simple polyhedron. Put
$$
\Omega^*_0(X,\alpha)_q=\omega_0^{-2\chi (X)+\chi (\partial X)}
\Omega_0(X,\alpha)_q,
$$
$$
\Omega^*_1(X,\alpha)_q=\omega^{-2\chi (X)+\chi (\partial X)}
\Omega_1(X,\alpha)_q,
$$
and
$$
\Omega^*_e(X,\alpha)_q=\omega^{-2\chi (X)+\chi (\partial X)}
(\Omega_0(X,\alpha)_q+\Omega_2(X,\alpha)_q).
$$

\th{Lemma}{2.}{Let $X\/$ be a simple 2-polyhedron and $\alpha\/$ be a
coloring of $\partial X$. Then the numbers $\Omega^*_0(X,\alpha)_q$,
$\Omega^*_1(X,\alpha)_q$, and
$\Omega^*_e(X,\alpha)_q$  are invariant under $\cal B$. }

\Proof It follows immediately from the definition of the number
 $|X,\varphi|_q$ that the number
$\Omega^*_0(X,\alpha)_q$ is invariant under $\cal B$ if
$$
\omega_0^2=w_j^{-2} \sum_{\scriptstyle k,l\equiv 0(2) \atop
\scriptstyle k,l: \{ j,k,l\}\in adm}
w_k^2 w_l^2
$$
for any even $j\in {\bf Z}_{r-1}$, and
$\Omega^*_1(X,\alpha)_q$, $\Omega^*_e(X,\alpha)_q$ are  invariant under
$\cal B$ if
$$
\omega^2=w_j^{-2} \sum_{k,l: \{ j,k,l\}\in adm }
w_k^2 w_l^2
$$
for any  $j\in {\bf Z}_{r-1}$.

The second equality is proved in [TV]. The proof of the first one is
 similar. First of all, let us check that
$$
w_j^{-2} \sum_{\scriptstyle k,l\equiv 0(2) \atop
\scriptstyle k,l: \{ j,k,l\}\in adm}
w_k^2 w_l^2=
w_0^{-2} \sum_{\scriptstyle s\equiv 0(2) \atop
\scriptstyle 0\le s \le r-2 }
w_s^4\leqno(*)
$$
for any even number $j\in {\bf Z}_{r-1}$.

Let $T\/$ be a polyhedron  obtained from a disk $D^2\/$ by attaching
one semidisk along a diameter of $D^2$. The polyhedron $T\/$ consists of
three 2-cells $\Gamma _1, \Gamma_2, \Gamma_3$.
Let a polyhedron $T_i$ is obtained from $T\/$ by attaching a 2-disk
along a circle that belongs to the 2-cell $\Gamma_i$, where $i=1\/$ or 2.

For any $j\in {\bf Z}_{r-1}$ we define a coloring $\beta$ of
$\partial T_1\/$ and $\partial T_2$ as  follows:
$\beta(\Gamma_1)=\beta(\Gamma_3)=j\/$, $\beta(\Gamma_2)=0$.

 By definition, we have
$$
\eqalign{
\Omega_0(T_1,\beta)_q&=w_0^2 w_j^2
\sum_{\scriptstyle k,l\equiv 0(2) \atop
\scriptstyle k,l: \{ j,k,l\}\in adm}
w_k^2 w_l^2,\cr
\Omega_0(T_2,\beta)_q&=w_j^4
\sum_{\scriptstyle s\equiv 0(2) \atop
\scriptstyle 0\le s \le r-2 }
w_s^4.\cr}
$$
Note that  $T_1$ and $T_2$ are connected by $\cal L$-move, therefore
$\Omega_0(T_1,\alpha)_q=\Omega_0(T_2,\alpha)_q$. This gives us the
equality  $(*)$.

Clearly, $w_0=1$. Thus we have to prove that
$$
\sum_{t=0}^{[r/2]-1} w_{2t}^4 =- r/(q-q^{-1})^2.
$$
The proof of this equality is straightforward.
\vskip4pt plus2pt

\th{Corollary}{2.}{Let $X$ be a simple spine of a 3-manifold $M$.
Then the numbers $TV^*_0(M)_q=\Omega^*_0(X)_q$,
$TV^*_1(M)_q=\Omega^*_1(X)_q$,   and
$TV^*_e(M)=\Omega^*_e(X)_q$ are invariants of $M$ under removing of
3-balls.}
\vskip4pt plus2pt

We can define the invariants $TV_0^*$, $TV_1^*$, and $TV_e^*$
  on
a triangulation of a 3-manifold $M\/$ like the Turaev--Viro invariants were
defined in [TV]. For simplicity we will restrict ourselves to the case
of closed 3-manifolds  only.

Let $M\/$ be a closed triangulated 3-manifold. Let $a\/$ be the number of
vertices of $M\/$, let $E_1,\ldots , E_b$ be the edges of $M\/$, and let
$T_1,\ldots , T_d$ be the 3-simplexes of $M\/$. By a {\it coloring\/}
of $M\/$ we
mean an arbitrary mapping
$\varphi\colon\{ E_1,\ldots , E_b\}\to {\bf Z}_{r-1}$.
A coloring $\varphi$ of $M$ is called {\it admissible} if for any 2-symplex
$A$ of $M$ the colors of the three edges of $A$ form an admissible triple.
Denote the set of admissible colorings of $M\/$ by $Adm(M)$.
We will denote by $Adm_0(M)$ the set of admissible colorings of $M\/$
by even numbers, and by $Adm_e(M)$ the set of admissible colorings of $M\/$
such that
$$
v-t+f\equiv 0 \pmod{2},
$$
where $v\/$ is the number of 3-simplexes  containing an edge
colored by an odd number,
$t\/$ is the number of 2-simplexes
containing an edge
colored by an odd number,
and $f\/$ is the number of  edges colored by odd numbers.
Note that $Adm_0(M)\subset Adm_e(M)$.
Set $Adm_1(M)=Adm(M)-Adm_e(M)$.
A 6-tuple
$\pmatrix{i&j&k\cr
          l&m&n\cr}$
is called a {\it color\/} of a 3-simplex $T_s$ if $i,j,k\/$ are the colors
of edges of some 2-face of $T_s$ and $(i,l)$, $(j,m)$, $(k,n)$ are the
pairs of  colors of opposite edges of $T_s$. Let
$$
|T_s^{\varphi}|_q=
\left\vert\matrix{i&j&k\cr
                  l&m&n\cr}\right\vert_q.
$$
For $\varphi\in Adm(M)$ put
$$
|M,\varphi|_q = \prod_{i=1}^b w^2_{\varphi (E_i),q}
\prod_{s=1}^d|T^{\varphi}_s|_q.
$$

\th{Preposition}{1.}{For any closed triangulated 3-manifold $M$ we have
$$
\eqalign{
TV^*_0(M)_q &= \omega_0^{-2a}
\sum_{\varphi\in Adm_0(M)} |M,\varphi|_q,\cr
TV^*_1(M)_q &= \omega^{-2a}
\sum_{\varphi\in Adm_1(M)} |M,\varphi|_q,\cr
TV^*_e(M)_q &= \omega^{-2a}
\sum_{\varphi\in Adm_e(M)} |M,\varphi|_q.\cr}
$$
}

\Proof
Let $X\/$ be the union of the closed barycentric stars of the edges of
$M\/$. It is obvious that $X\/$ is a special polyhedron. By a finite
number of ${\cal M}^{\pm 1}, {\cal L}^{\pm 1}$ and ${\cal B}^{- 1}$
moves on $X\/$ we get a simple spine of $M$. Each coloring $\varphi\/$
of $M\/$ induces a dual coloring $\varphi^*\/$ of $X\/$, and it is
easy to check that $|M,\varphi|_q=|X,\varphi^*|_q$ and $\chi (X)=a$,
which establishes the formulas.
\vskip4pt plus2pt

\section{6. The values of the summand-invariants.}
Here we express the numbers \break $TV_0(M)_q+TV_2(M)_q$ and $TV_1(M)_q$ via
the Turaev--Viro invariants.

Let $X\/$ be a special polyhedron. Fix a number $r\ge 3$ and a coloring
$\varphi\in Adm(X)$. A vertex $v\/$ of the colored polyhedron $X\/$ is
called a {\it switch-vertex\/} if the sum of all odd numbers in the color
of $v\/$ is congruent to 2 modulo 4.

\th{Lemma}{3.}
{ Let $X$ be a special polyhedron. Then for any
$\varphi\in Adm(X)$ we have
$$
|X,\varphi|_q=(-1)^{\chi (S(\varphi))+x} |X,\varphi|_{-q},
$$
where $x$ is the number of the switch-vertices of $X$.
}

\Proof
It is easy to see that
$$[n]_q=(-1)^{n-1}[n]_{-q},$$
$$[n]_q! = (-1)^{n(n-1)/2}[n]_{-q}!, $$
$$w^2_{i,q}=(-1)^i w^2_{i, -q},$$

$$
\Delta^2_q (i,j,k)=
\cases{
\Delta^2_{-q}(i,j,k), & if $i,j,k$ are even,\cr
-\Delta^2_{-q}(i,j,k), &  otherwise.\cr}
$$

Let      $\pmatrix{i_1&i_2&i_3\cr
             i_4&i_5&i_6\cr}$
be a color of a vertex $v$ of $X$ under $\varphi$.
Then we have
$$
(**)\qquad \left[\matrix{i_1&i_2&i_3\cr
                        i_4&i_5&i_6\cr}\right]_q =
(-1)^{{{1}\over {2}} \sum_{\scriptstyle s,t=1\atop
                        \scriptstyle s\le t}^6 i_s i_t}
\left[\matrix{i_1&i_2&i_3\cr
              i_4&i_5&i_6\cr}\right]_{-q}.
$$

There are three possibilities for the color of $v$.

1) Each number in the color of $v$ is even ({\it even vertex}). Then the
 sign in $(**)$ is plus.

2) There are four odd numbers in the color of $v$ ({\it fourfold vertex}).
Let $i_1, i_2, i_4, i_5$
be the odd numbers, then from $(**)$ we have
$$
 \left[\matrix{i_1&i_2&i_3\cr
                        i_4&i_5&i_6\cr}\right]_q =
(-1)^{{{i_1+i_2+i_4+i_5}\over {2}}+1}
\left[\matrix{i_1&i_2&i_3\cr
              i_4&i_5&i_6\cr}\right]_{-q}.
$$

Hence if $v$ is a switch-vertex, then the sign in $(**)$ is plus,
otherwise minus.

3) There are three odd  and three even numbers in the color of $v$
({\it threefold vertex}). Let
$i_1, i_2, i_3$ be the even numbers, then from $(**)$ we have
$$
 \left[\matrix{i_1&i_2&i_3\cr
                        i_4&i_5&i_6\cr}\right]_q =
(-1)^{{{i_1+i_2+i_3}\over {2}}+1}
\left[\matrix{i_1&i_2&i_3\cr
              i_4&i_5&i_6\cr}\right]_{-q}.
$$

By a {\it cost} of an edge $E$ we mean the half-sum of the numbers from
a color of $E$. Let $E'$ be a half-edge of an edge $E$. By a {\it cost}
of the half-edge $E'$ we mean the cost of  $E$.
We call an edge or a half-edge  {\it bad} if its color is
{\it even} (that is all three numbers in the color of the edge are even)
and its cost
is even. Let us show that the number of threefold vertices with a bad
half-edge is even.
It is sufficient to prove that the number of the bad half-edges
incident to an even
vertex is even, but this statement follows from the fact that the sum
of all 4 costs of the half-edges  incident to an even vertex is even.
Hence we can think that for any threefold vertex the sign in $(**)$ is plus.

Let us denote the number of the odd-colored edges of $X$ by $e$,
the number of the threefold
vertices by $n_3$, and the number of the fourfold vertices by $n_4$.
Denote the number of the odd colored 2-components of $X$ by $c$.
Then we have $c=\chi(S(\varphi))-n_3-n_4+e$. Hence
$$
|X,\varphi|_q =(-1)^{c+n_4-x+e}
|X,\varphi|_{-q}=(-1)^{\chi(S(\varphi))-n_3-x}|X,\varphi|_{-q}.
$$
It is easy to see that for any admissible coloring $\varphi$ of a special
polyhedron $X$ the number $n_3$ is even. This finishes the proof.
\vskip4pt plus2pt

Let $SX\/$ be the set of  singular points of $X\/$. Note that $SX\/$ is
a regular graph of degree~4. Denote by $V\/$ the set of vertices of $X\/$,
by $N(V,SX)\/$ a closed regular neighborhood of $V\/$ in $SX\/$,
and by $N(V,X)\/$
a closed regular neighborhood of $V\/$ in $X\/$.
The intersection of the union
of the open edges with each connected component of $N(V,SX)\/$ consists
of 4 half-open 1-cells, which are called {\it thorns\/}. The intersection
of the union of the open 2-cells with each connected component of $N(V,X)$
consists of six half-open 2-cells, which are called {\it wings\/}.

Let $v$ be a vertex of $X$, and let $N(v,M)$ be a closed regular
neighborhood of $v$. Chose a thorn $t$ in $N(v,M)$ and a small normal disk
$D$ for it. Any orientation $\alpha$ of $N(v,M)$ induces an orientation
$\alpha |_D$ of $D$ according to the following convention: $\alpha |_D$
together with the outward orientation of $t$ should give the orientation
$\alpha$. Note that $\alpha |_D$ induces a cyclic order on the set of
wings adjacent to $t$.

Regular neighborhood $N(V,M)$ consists of 3-balls $N(v,M),\quad v\in V$.
Choose orientations for the 3-balls.
Let $E\/$ be an edge of $X\/$. It contains two thorns $t_1, t_2$.
Let $W_1^{(i)}, W_2^{(i)}, W_3^{(i)}\/$ be the wings adjacent to $t_i$,
where $i=1,2$.
As above, the orientation of $N(V,M)$ induces a cyclic order on the set
$$
\{W_1^{(i)}, W_2^{(i)}, W_3^{(i)}\},\qquad {for}\quad i=1,2.
$$
 The 2-cells of $X\/$ determine the
natural bijection
$$
f\colon\{ W_1^{(1)},W_2^{(1)},W_3^{(1)}\}\to
\{ W_1^{(2)},W_2^{(2)},W_3^{(2)}\}.
$$
 We shall say that the edge $E\/$ is
{\it odd\/} if the bijection $f\/$ preserves the cyclic order on the wings,
and {\it even\/} otherwise.

\th{Theorem}{2.}{Let $X$ be a special spine of a 3-manifold $M$.
Then for any $\varphi\in Adm(X)$ we have
$$
|X,\varphi|_q=(-1)^{\chi (S(\varphi))} |X,\varphi|_{-q}.
$$}

\Proof
Let $x\/$ be the number of switch-vertices of the pair $(X,\varphi)$.
By lemma 3 it is sufficient to prove that this number is even.

Consider the coloring $\overline\varphi\colon\{ \Gamma_1,\ldots ,\Gamma _b
\} \to {\bf Z}_4$ such that
$$
\overline\varphi (\Gamma _i)=
\cases{0,& if $\varphi (\Gamma_i)\equiv 0 \pmod{2}$,\cr
       1,& if $\varphi (\Gamma_i)\equiv 1 \pmod{4}$,\cr
       3,& if $\varphi (\Gamma_i)\equiv 3 \pmod{4}$,\cr}
$$
for any $1\le i\le b$.

Fix an orientation of $N(V,M)$. Then each edge of $SX\/$ becomes odd or
even. Let $G\/$ be the union of the edges of $X\/$ with the colors
$\{ 0,1,3 \}$ under the coloring $\overline\varphi$.
 Let
$\Omega_1,\ldots ,\Omega _p$ be the middle points of the odd edges of $G\/$.
Consider a graph
$G'$. The set of vertices of $G'\/$ consists of the vertices of $G\/$ and
of the  points $\Omega_1,\ldots ,\Omega _p$.
The set of edges of $G'\/$ consists of the even edges of $G\/$ and of the halves
of the odd edges of $G\/$. So each odd edge of $G\/$ gives 2 edges in $G'$.
The orientation of $N(V,M)$ and the coloring $\overline\varphi$ give
the orientation of the graph $G'$. Let $v_1,\ldots ,v_t$ be the vertices of
$G'$. We will denote by $a_i$ the number of  incoming and by $b_i$
the number of outgoing edges for the vertex $v_i$. We have
$(a_i-b_i)\equiv 2\pmod{4}$ iff $v_i$ is either the switch-vertex or the
middle point of an odd edge, and $(a_i-b_i)\equiv 0\pmod{4}$ otherwise.
The number
of vertices with the condition $(a_i-b_i)\equiv 2\pmod{4}$  is even for any
oriented graph, because $\sum_{i=1}^t (a_i-b_i)=0$.

It remains to prove that the number $p\/$ is even. Let $\theta\/$ be the
number of odd edges of $X\/$ with the color $\{ 0,1,1\}$ under the coloring
$\overline\varphi\/$. Then 1-colored (by $\overline\varphi\/$) 2-cells
pass  $(2\theta +p)$ times along the odd edges of $X\/$.
Note that each 2-component of $X\/$ passes along the odd edges of $X\/$
even number of times (this is true for every  special spine;
see, for instance, [F]).
Therefore the number
$(2\theta +p)$ is even and $p\/$ is even.
\vskip4pt plus2pt

\remar{Remark\ {4.}\ }{In the case of an orientable 3-manifold this theorem was
proved in [S1].}
\vskip4pt plus2pt

\th{Corollary}{3.}{For any 3-manifold $M$ and any  $q$ we have
$$
\eqalign{
TV_N(M)_q &=(-1)^N TV_N(M)_{-q},\quad {where}\quad
N\in\{ 0,1,2\},\cr
TV_0(M)_q &+TV_2(M)_q={1\over 2} (TV(M)_q+TV(M)_{-q}),\cr
TV_1(M)_q &={1\over 2} (TV(M)_q-TV(M)_{-q}).\cr}
$$
}

\remar{Remark\ {5.}\ }{In the papers [S1] and [S2] we used the parameter
$-\overline{q}$ instead of $-q$, but it is easy to see that
$[n]_q=[n]_{\overline{q}}$.
}
\vskip4pt plus2pt

\section{7. The tables.}
Below we present the summand-invariants $TV_N(M)_q$ and the Turaev--Viro
invariants $TV^*(M)_q$ with $r\le 7$ for the manifolds $S^3, {\bf R}P^3,
L_{3,1}, L_{4,1}, L_{5,1}, L_{5,2}, L_{6,1}$, $L_{7,2}, L_{8,3}, L_{9,2},
L_{10,3}, L_{11,4}, L_{12,5}, L_{13,5}, S^3/Q_8, S^3/Q_{12}$, where
$S^3/G\/$ denotes the quotient space of the sphere $S^3\/$ by a linear free
action of a finite nonabelian group $G$. These are all closed irreducible
 orientable 3-manifolds, having a special spine with $\le 3$ vertices.

Each summand invariant is presented by a polynomial on $q\/$ (here $q$ is a
primitive root of unity of degree $2r$) with integer
coefficients, and by
evaluation of the polynomial  at $q=e^{i\pi /r}$. Note that each
coefficient in the polynomial is a separate invariant.
The invariants from the tables are
related by the equality
$$
TV^*(M)_q={{(q-q^{-1})^2}\over {2r}}(TV_0(M)_q+TV_1(M)_q+TV_2(M)_q).
$$

\nine
\bigskip
Table 1: Invariants for $S^3$

\medskip
\settabs\+\indent & r\qquad &00000000000&**********&1111111111&**********&22222222222&**********&\cr
\+&\negthinspace \it r&$TV_0(M)_q$&&$TV_1(M)_q$&&$TV_2(M)_q$&&$TV^*(M)_q$\cr
\smallskip
\hrule height0.4pt width30truecc
\smallskip
\+&3&1&=1.000 &0&=0.000 &0&=0.000 &0.500\cr
\+&4&1&=1.000 &0&=0.000 &0&=0.000 &0.250\cr
\+&5&1&=1.000 &0&=0.000 &0&=0.000 &0.138\cr
\+&6&1&=1.000 &0&=0.000 &0&=0.000 &0.083\cr
\+&7&1&=1.000 &0&=0.000 &0&=0.000 &0.054\cr
\bigskip

\break

\bigskip
Table 2: Invariants for ${\bf R}P^3$
\medskip
\settabs\+\indent &r\quad &$-2q^5+q^4-q^3+2q^2+3$&********&$2q^5+q^4-q^3+2q^2+3$&*********&2&********&\cr
\+&\negthinspace \it r&$TV_0(M)_q$&&$TV_1(M)_q$&&$TV_2(M)_q$&&$TV^*(M)_q$\cr
\smallskip
\hrule height0.4pt width30truecc
\smallskip
\+&3&1&=1.000 &$-1$&=$-1.000$ &0&=0.000 &0.000\cr
\+&4&2&=2.000 &$q^3-q$&=$-1.414$ &0&=0.000 &0.146\cr
\+&5&$-q^3+q^2+2$&=2.618 &$q^3-q^2-2$&=$-2.618$  &0&=0.000 &0.000\cr
\+&6&4&=4.000 &$2q^3-4q$ &=$-3.464$ &0&=0.000 &0.045\cr
\+&7 &$-2q^5+q^4-q^3+2q^2+3$ &=5.049 &$2q^5-q^4+q^3-2q^2-3$ &=$-5.049$ & 0 &=0.000  &0.000\cr
\bigskip

\bigskip
Table 3: Invariants for $L_{3,1}$

\bigskip
\settabs\+\indent & r\qquad &0000000000000&**********&1111111111&**********&22222222222&**********&\cr
\+&\negthinspace \it r&$TV_0(M)_q$&&$TV_1(M)_q$&&$TV_2(M)_q$&&$TV^*(M)_q$\cr
\smallskip
\hrule height0.4pt width30truecc
\smallskip
\+&3&1&=1.000 &0&=0.000 &0&=0.000 &0.500\cr
\+&4&1&=1.000 &0&=0.000 &0&=0.000 &0.250\cr
\+&5&$-q^3+q^2+2$&=2.618 &0&=0.000 &0&=0.000 &0.362\cr
\+&6&3&=3.000 &0&=0.000 &0&=0.000 &0.250\cr
\+&7&$-q^5+2q^2+2$&=3.247 &0&=0.000 &0&=0.000 &0.175\cr
\bigskip

\bigskip
Table 4: Invariants for $L_{4,1}$
\medskip
\settabs\+\indent & r\qquad &$-q^5+2q^2+2$ &*********&$TV_0(M)_q$&*********&$-q^5+2q^2+2$ &*********&\cr
\+&\negthinspace \it r&$TV_0(M)_q$&&$TV_1(M)_q$&&$TV_2(M)_q$&&$TV^*(M)_q$\cr
\smallskip
\hrule height0.4pt width30truecc
\smallskip
\+&3&1&=1.000 &0&=0.000 &1&=1.000 &1.000\cr
\+&4&2&=2.000 &0&=0.000 &0&=0.000 &0.500\cr
\+&5&1&=1.000 &0&=0.000 &1&=1.000 &0.276\cr
\+&6&4&=4.000 &0&=0.000 &0&=0.000 &0.333\cr
\+&7&$-q^5+2q^2+2$&=3.247 &0&=0.000 &$-q^5+2q^2+2$&=3.247  &0.349\cr
\bigskip

\bigskip
Table 5: Invariants for $L_{5,1}$
\medskip
\settabs\+\indent & r\qquad &$-2q^5+q^4-q^3+2q^2+3$ &*********&$2q^2+3$&*********&22222 &*********&\cr
\+&\negthinspace \it r&$TV_0(M)_q$&&$TV_1(M)_q$&&$TV_2(M)_q$&&$TV^*(M)_q$\cr
\smallskip
\hrule height0.4pt width30truecc
\smallskip
\+&3&1&=1.000 &0&=0.000 &0&=0.000 &0.500\cr
\+&4&1&=1.000 &0 &=0.000 &0&=0.000 &0.250\cr
\+&5&$-q^3+q^2+3$&=3.618 &0 &=0.000  &0&=0.000 &0.500\cr
\+&6&1&=1.000 &0 &=0.000 &0&=0.000 &0.083\cr
\+&7 &$-2q^5+q^4-q^3+2q^2+3$ &=5.049 &0 &=0.000 &0 &=0.000  &0.272\cr
\bigskip

\bigskip
Table 6: Invariants for $L_{5,2}$

\medskip
\settabs\+\indent & r\qquad &$-2q^5+q^4-q^3+2q^2+3$ &*********&$TV_0(M)_q$&*********&$TV_0(M)_q$&*********&\cr
\+&\negthinspace \it r&$TV_0(M)_q$&&$TV_1(M)_q$&&$TV_2(M)_q$&&$TV^*(M)_q$\cr
\smallskip
\hrule height0.4pt width30truecc
\smallskip
\+&3&1&=1.000 &0&=0.000 &0&=0.000 &0.500\cr
\+&4&1&=1.000 &0&=0.000 &0&=0.000 &0.250\cr
\+&5&0&=0.000 &0&=0.000 &0&=0.000 &0.000\cr
\+&6&1&=1.000 &0&=0.000 &0&=0.000 &0.083\cr
\+&7&$-2q^5+q^4-q^3+2q^2+3$&=5.049 &0&=0.000 &0&=0.000 &0.272\cr
\bigskip

\break

\bigskip
Table 7: Invariants for $L_{6,1}$
\medskip
\settabs\+\indent & r\qquad &$-2q^5+7$ &*********&$-q^3+q$&*********&$-6q^2+21$&*********&\cr
\+&\negthinspace \it r&$TV_0(M)_q$&&$TV_1(M)_q$&&$TV_2(M)_q$&&$TV^*(M)_q$\cr
\smallskip
\hrule height0.4pt width30truecc
\smallskip
\+&3&1&=1.000 &$-1$&=$-1.000$ &0&=0.000 &0.000\cr
\+&4&2&=2.000 &$-q^3+q$ &=1.414 &0&=0.000 &0.853\cr
\+&5&1&=1.000 &$-1$&=$-1.000$  &0&=0.000 &0.000\cr
\+&6&6&=6.000 &0 &=0.000 &0&=0.000 &0.500\cr
\+&7 &1  &=1.000 &$-1$&=$-1.000$ &0 &=0.000  &0.000\cr
\bigskip

\bigskip
Table 8: Invariants for $L_{7,2}$
\medskip
\settabs\+\indent & r\qquad &$-q^3+q^2+2$ &*********&$2q^2+3$&*********&22222 &*********&\cr
\+&\negthinspace \it r&$TV_0(M)_q$&&$TV_1(M)_q$&&$TV_2(M)_q$&&$TV^*(M)_q$\cr
\smallskip
\hrule height0.4pt width30truecc
\smallskip
\+&3&1&=1.000 &0&=0.000 &0&=0.000 &0.500\cr
\+&4&1&=1.000 &0 &=0.000 &0&=0.000 &0.250\cr
\+&5&$-q^3+q^2+2$&=2.618 &0 &=0.000  &0&=0.000 &0.362\cr
\+&6&1&=1.000 &0 &=0.000 &0&=0.000 &0.083\cr
\+&7 &0  &=0.000 &0 &=0.000 &0 &=0.000  &0.000\cr
\bigskip

\bigskip
Table 9: Invariants for $L_{8,3}$
\medskip
\settabs\+\indent & r\qquad &$-q^3+q^2+2$ &*********&$-q^3+2$&*********&$-q^3+q^2+2$&*********&\cr
\+&\negthinspace \it r&$TV_0(M)_q$&&$TV_1(M)_q$&&$TV_2(M)_q$&&$TV^*(M)_q$\cr
\smallskip
\hrule height0.4pt width30truecc
\smallskip
\+&3&1&=1.000 &0&=0.000 &1&=1.000 &1.000\cr
\+&4&2&=2.000 &0 &=0.000 &2&=2.000 &1.000\cr
\+&5&$-q^3+q^2+2$&=2.618 &0 &=0.000  &$-q^3+q^2+2$&=2.618 &0.724\cr
\+&6&4&=4.000 &0 &=0.000 &0&=0.000 &0.333\cr
\+&7 &1  &=1.000 &0 &=0.000 &1 &=1.000  &0.108\cr
\bigskip

\bigskip
Table 10: Invariants for $L_{9,2}$
\medskip
\settabs\+\indent & r\qquad &$-2q^5+q^4-q^3+2q^2+3$ &*********&$-q^3+q$&*********&$-6q^2+21$&*********&\cr
\+&\negthinspace \it r&$TV_0(M)_q$&&$TV_1(M)_q$&&$TV_2(M)_q$&&$TV^*(M)_q$\cr
\smallskip
\hrule height0.4pt width30truecc
\smallskip
\+&3&1&=1.000 &0 &=0.000 &0&=0.000 &0.500\cr
\+&4&1&=1.000 &0 &=0.000 &0&=0.000 &0.250\cr
\+&5&1&=1.000 &0 &=0.000  &0&=0.000 &0.138\cr
\+&6&3&=3.000 &0 &=0.000 &0&=0.000 &0.250\cr
\+&7 &$-2q^5+q^4-q^3+2q^2+3$ &=5.049 &0 &=0.000 &0 &=0.000  &0.272\cr
\bigskip

\bigskip
Table 11: Invariants for $L_{10,3}$
\medskip
\settabs\+\indent & r\qquad &$-q^5+q^2+2$ &*********&$-q^5+q^2+2$&*********&$-6q8821$&*********&\cr
\+&\negthinspace \it r&$TV_0(M)_q$&&$TV_1(M)_q$&&$TV_2(M)_q$&&$TV^*(M)_q$\cr
\smallskip
\hrule height0.4pt width30truecc
\smallskip
\+&3&1&=1.000 &$-1$&=$-1.000$ &0&=0.000 &0.000\cr
\+&4&2&=2.000 &$-q^3+q$ &=1.414 &0&=0.000 &0.853\cr
\+&5&0&=0.000 &0 &=0.000  &0&=0.000 &0.000\cr
\+&6&4&=4.000 &$-2q^3+4q$ &=3.464 &0  &=0.000 &0.622\cr
\+&7 &$-q^5+q^2+2$  &=3.247 &$q^5-q^2-2$  &=$-3.247$ &0 &=0.000  &0.000\cr
\bigskip

\break

\bigskip
Table 12: Invariants for $L_{11,4}$
\medskip
\settabs\+\indent & r\qquad &$-q^5+q^2+2$ &*********&$-q^5+q^2+2$&*********&$-6q8821$&*********&\cr
\+&\negthinspace \it r&$TV_0(M)_q$&&$TV_1(M)_q$&&$TV_2(M)_q$&&$TV^*(M)_q$\cr
\smallskip
\hrule height0.4pt width30truecc
\smallskip
\+&3&1&=1.000 &0  &= 0.000 &0&=0.000 &0.500\cr
\+&4&1&=1.000 &0 &=0.000 &0&=0.000 &0.250\cr
\+&5&1&=1.000 &0 &=0.000  &0&=0.000 &0.138\cr
\+&6&1&=1.000 &0 &=0.000 &0  &=0.000 &0.083\cr
\+&7 &$-q^5+q^2+2$  &=3.247 &0  &=0.000 &0 &=0.000  &0.175\cr
\bigskip

\bigskip
Table 13: Invariants for $L_{12,5}$
\medskip
\settabs\+\indent &r\quad &$-2q^5+q^4-q^3+2q^2+3$&********&2&********&$-2q^5+q^4-q^3+2q^2+3$&********&\cr
\+&\negthinspace \it r&$TV_0(M)_q$&&$TV_1(M)_q$&&$TV_2(M)_q$&&$TV^*(M)_q$\cr
\smallskip
\hrule height0.4pt width30truecc
\smallskip
\+&3&1&=1.000 &0 &=0.000 &1&=1.000 & 1.000\cr
\+&4&2&=2.000 &0 &=0.000 &0&=0.000 & 0.500\cr
\+&5&$-q^3+q^2+2$ &=2.618 &0 &=0.000  &$-q^3+q^2+2$ &=2.618 & 0.724\cr
\+&6&6&=6.000 &0 &=0.000 &6&=6.000 & 1.000\cr
\+&7 &$-2q^5+q^4-q^3+2q^2+3$ &=5.049 &0 &=0.000 &$-2q^5+q^4-q^3+2q^2+3$ &=5.049  & 0.543\cr
\bigskip

\bigskip
Table 14: Invariants for $L_{13,5}$

\medskip
\settabs\+\indent & r\qquad &00000000000&**********&1111111111&**********&22222222222&**********&\cr
\+&\negthinspace \it r&$TV_0(M)_q$&&$TV_1(M)_q$&&$TV_2(M)_q$&&$TV^*(M)_q$\cr
\smallskip
\hrule height0.4pt width30truecc
\smallskip
\+&3&1&=1.000 &0&=0.000 &0&=0.000 &0.500\cr
\+&4&1&=1.000 &0&=0.000 &0&=0.000 &0.250\cr
\+&5&$-q^3+q^2+2$&=2.618 &0&=0.000 &0&=0.000 &0.362\cr
\+&6&1&=1.000 &0&=0.000 &0&=0.000 &0.083\cr
\+&7&1&=1.000 &0&=0.000 &0&=0.000 &0.054\cr
\bigskip

\bigskip
Table 15: Invariants for $S^3/Q_8$
\medskip
\settabs\+\indent & r\qquad &$-2q^5+2q^2+7$ &*********&$-q^3+2$&*********&$-6q^5+6q^2+21$&*********&\cr
\+&\negthinspace \it r&$TV_0(M)_q$&&$TV_1(M)_q$&&$TV_2(M)_q$&&$TV^*(M)_q$\cr
\smallskip
\hrule height0.4pt width30truecc
\smallskip
\+&3&1&=1.000 &0&=0.000 &3&=3.000 &2.000\cr
\+&4&4&=4.000 &0 &=0.000 &6&=6.000 &2.500\cr
\+&5&$-q^3+q^2+4$&=4.618 &0 &=0.000  &$-q^3+3q^2+12$&=13.854 &2.553\cr
\+&6&10&=10.000 &0 &=0.000 &18&=18.000 &2.333\cr
\+&7 &$-2q^5+2q^2+7$  &=9.494 &0 &=0.000 &$-6q^5+6q^2+21$ &=28.482  &2.043\cr
\bigskip

\bigskip
Table 16: Invariants for $S^3/Q_{12}$
\medskip
\settabs\+\indent &r\quad &$-2q^5+q^4-q^3+2q^2+3$&********&2&********&$-2q^5+q^4-q^3+2q^2+3$&********&\cr
\+&\negthinspace \it r&$TV_0(M)_q$&&$TV_1(M)_q$&&$TV_2(M)_q$&&$TV^*(M)_q$\cr
\smallskip
\hrule height0.4pt width30truecc
\smallskip
\+&3&1&=1.000 &0 &=0.000 &1&=1.000 & 1.000\cr
\+&4&2&=2.000 &0 &=0.000 &0&=0.000 & 0.500\cr
\+&5&$-q^3+q^2+4$ &=4.618 &0 &=0.000  &$-q^3+q^2+4$ &=4.618 & 1.276\cr
\+&6&10&=10.000 &0 &=0.000 &6&=6.000 & 1.333\cr
\+&7 &$-2q^5+q^4-q^3+2q^2+5$ &=7.049 &0 &=0.000 &$-2q^5+q^4-q^3+2q^2+5$ &=7.049  & 0.758\cr
\bigskip

\break

\references{BHMV}{
\item{[BHMV]} C. \spa{Blanchet}, N. \spa{Habegger}, G. \spa{Masbaum},
P. \spa{Vogel},
{\it Remarks on the three-manifold invariants $\theta_p$}\/,
'Operator Algebras, Mathematical Physics, and Low Dimensional Topology'
(Nato Workshop July 1991), Vol.~5, 39--59.

\item{[Ca]} B.~G. \spa{Casler},
{\it An imbedding theorem for connected 3-manifolds with boundary}\/,
Math. Soc. 16 (1965), 559--566.

\item{[F]} A.~T. \spa{Fomenko},
{\it Visual geometry and topology}\/,
Moscow University, 1992 (in Russian).

\item{[KL1]} L.~H. \spa{Kauffman} and S. \spa{Lins},
{\it Computing Turaev--Viro invariants for 3-manifolds}\/,
Manuscripta Math. 72 (1991), 81--94.

\item{[KL2]} L.~H. \spa{Kauffman} and S. \spa{Lins},
{\it Temperley--Lieb recoupling theory and invariants of 3-manifolds}\/,
Princeton University Press, Princeton, N.~J., 1994.

\item{[KM]} R.~C. \spa{Kirby}, P. \spa{Melvin},
{\it The 3-manifold invariants of Witten and Reshetikhin--Turaev for
$sl(2,{\bf C})$}\/,
Inv. Math. 105 (1991), 473--545.

\item{[M]} S.~V. \spa{Matveev},
{\it Transformations of special spines and the Zeeman conjecture}\/,
Math. USSR Izvestia, Vol.~31, 2, 1988, 423--434.

\item{[N]} J.~R. \spa{Neil},
{\it Combinatorial calculation of the various normalizations of the
Witten invariants for 3-manifolds}\/,
J. of Knot Theory and Its Ramifications, Vol.~1, 4 (1992),
407--449.

\item{[P]} R. \spa{Piergallini},
{\it Standard moves for standard polyhedra and spines}\/,
Rend. Circ. Mat. Palermo, Vol.~37, 18 (1988), 391--414.

\item{[S1]}  M. \spa{Sokolov},
{\it On Turaev--Viro invariants for 3-manifolds}\/,
VINITI preprint 583--B93 (in Russian).

\item{[S2]}  M. \spa{Sokolov},
{\it The Turaev--Viro invariant for 3-manifolds is a sum of three
invariants}\/,
Can. Math. Bull., (to appear).

\item{[T1]} V.~G. \spa{Turaev},
{\it Topology of shadows}\/.
Preprint, 1991.

\item{[T2]} V.~G. \spa{Turaev},
{\it Quantum invariants of knots and 3-manifolds}\/,
Walter de Gruyter, Berlin, New York, 1994.

\item{[TV]} V.~G. \spa{Turaev} and O.~Y. \spa{Viro},
{\it State sum invariants of 3-manifolds and quantum 6j-symbols}\/,
Topology, 31 (1992), 866--902.

\item{[W]} K. \spa{Walker},
{\it On Witten's 3-manifold invariants}\/,
Preprint, 1991.}

\bye